\documentclass[12pt]{amsart}
\usepackage{amssymb,latexsym,comment,url}

\newcommand{\Char}{\operatorname{char}}

\newcommand{\Gal}{\operatorname{Gal}}
\newcommand{\GL}{\operatorname{GL}}

\newcommand{\CC}{\operatorname{\mathcal C}}
\newcommand{\PP}{{\mathbb P}}
\newcommand{\Q}{{\mathbb Q}}

\newcommand{\Z}{{\mathbb Z}}

\newenvironment{Proof}{\par\noindent{\sc Proof:}}%
                      {\hspace*{\fill}\nobreak$\Box$\par\medskip}
                       {\hspace*{\fill}\nobreak$\Box$\par\medskip}

\newtheorem{Proposition}{Proposition}[section]
\newtheorem{Theorem}[Proposition]{Theorem}
\newtheorem{Lemma}[Proposition]{Lemma}
\newtheorem{Corollary}[Proposition]{Corollary}

\theoremstyle{definition}
\newtheorem{Definition}[Proposition]{Definition}
\newtheorem{Remark}[Proposition]{Remark}

\renewcommand{\baselinestretch}{1.1}

\begin{document}

\title[Regular models of quadratic twists]%
{Minimal regular models of quadratic twists of genus two curves}

\author{Mohammad~Sadek}
\address{American University in Cairo, Mathematics and Actuarial Science Department, AUC Avenue, New Cairo, Egypt}
\email{mmsadek@aucegypt.edu}
\let\thefootnote\relax
\let\thefootnote\relax

\begin{abstract}
Let $K$ be a complete discrete valuation field with ring of integers $R$ and residue field $k$ of characteristic $p>2$. We assume moreover that $k$ is algebraically closed. Let $C$ be a smooth projective geometrically connected curve of genus $2$. If $K(\sqrt{D})/K$ is a quadratic field extension of $K$ with associated character $\chi$, then $C^{\chi}$ will denote the quadratic twist of $C$ by $\chi$. Given the minimal regular model $\mathcal X$ of $C$ over $R$, we determine the minimal regular model of the quadratic twist $C^{\chi}$. This is accomplished by obtaining the stable model $\CC^{\chi}$ of $C^{\chi}$ from the stable model $\CC$ of $C$ via analyzing the Igusa and affine invariants of the curves $C$ and $C^{\chi}$, and calculating the degrees of singularity of the singular points of $\CC^{\chi}$.

\end{abstract}

\maketitle


\section{Introduction}
 Let $R$ be a complete discrete valuation ring, $K$ its field of fractions, and $k$ its residue field. We assume moreover that $k$ is algebraically closed and $\Char k\ne2$. Let $\nu:K\to\Z\cup\{\infty\}$ be the normalized valuation of $K$.

  Let $C$ be a smooth projective geometrically connected curve of genus $g>0$ over $K$. There exists a regular $R$-curve $X$ whose generic fiber is isomorphic to $C$ and such that if $X'$ is another regular $R$-curve whose generic fiber is isomorphic to $C$ then $X'$ dominates $X$. $X$ is said to be the minimal regular model of $C$. The type of the minimal regular model is determined by the structure of its special fiber $X_k$ which we are going to call the reduction type of $X$.

    If $C$ is an elliptic curve, then given a Weierstrass equation describing $C$, Tate's algorithm produces the minimal regular model of $C$ by analyzing the $a$- and $b$-invariants of $C$, see \cite[Chapter IV, \S 9]{Silverman}. If $C$ has genus $2$, then the complete classification of $X_k$ can be found in \cite{NamikawaUeno}. The reduction type of $X$ takes one out of more than $120$ possibilities.

   Let $C$ be a curve of genus $2$. Let $L/K$ be the smallest field extension over which $C$ admits a stable model. In \cite{Liumodelesminimaux}, using the stable model of $C$ together with a careful analysis of the Igusa and affine invariants attached to $C$, Liu reproduced the reduction type of the minimal regular model of $C$.

If $K(\sqrt{D})$ is a quadratic extension of $K$ with associated character $\chi$, we denote the corresponding quadratic twist of the curve $C$ by $C^{\chi}$. If $C$ is an elliptic curve, given the reduction type of the minimal regular model $X$ of $C$, one can find the reduction type of the minimal regular model $X^{\chi}$ of $C^{\chi}$. Indeed, if $\nu(D)$ is even, then the reduction type of $X^{\chi}$ is the same as the one of $X$. If $\nu(D)$ is odd, then the complete list of reduction types of $X^{\chi}$ can be found in \cite{Comalada}. For example, if $\Char k\ne 2$, then

\begin{center}
\begin{tabular}{|c|c|c|c|c|c|c|c|c|c|c|}
  \hline
  type of $X_k$ & ${\rm I}_0$ & ${\rm I}_n$ & ${\rm II}$ & ${\rm III}$ & ${\rm IV}$ & ${\rm I}_0^*$ & ${\rm I}_n^*$ & ${\rm IV}^*$ &${\rm III}^*$&${\rm II}^*$\\
  \hline
  type of $X^{\chi}_k$ & ${\rm I}_0^*$ & ${\rm I}_n^*$ & ${\rm IV}^*$ & ${\rm III}^*$ & ${\rm II}^*$ & ${\rm I}_0$ & ${\rm I}_n$ & ${\rm II}$&${\rm III}$&${\rm IV}$ \\
  \hline
\end{tabular}
\end{center}

If $C$ has genus $g>1$, then one may pose the same question on the reduction type of $X^{\chi}$. In \cite{Sadek}, if $C$ has genus $g>0$, the description of the special fiber of $X^{\chi}$ is given if $X$ is smooth. In this article, given the reduction type of $X$, we display the reduction type of $X^{\chi}$ if the genus of $C$ is $g=2$.
We first investigate the Igusa and affine invariants attached to $C^{\chi}$ and link these invariants to the ones attached to $C$. This enables us to describe the stable model $\CC^{\chi}$ of $C^{\chi}$ and compute the degrees of singularity of its ordinary double points. Then one will be able to construct the minimal regular model $X^{\chi}$ from $\CC^{\chi}$.

\section{Invariants of genus two curves}
Let $C$ be a smooth projective geometrically connected genus two curve defined over a field $K$. If $\Char K\ne 2$, then $C$ is defined by a hyperelliptic equation of the form $y^2=P(x)$ where $P(x)$ is a polynomial in $K[x]$ of degree $5$ or $6$ with no repeated roots. If $K$ is the fraction field of some ring $R$, one may assume that $P(x)\in R[x]$. Let $f:C\to\PP_K^1$ be a finite separable morphism of degree $2$. Let $\sigma$ be the hyperelliptic involution of $C$. A point $x\in C(K)$ is a ramification point of $f$ if $\sigma(x)=x$. In particular, the ramification points of $f$ are the zeros of $P(x)$, plus the point at infinity if $\deg P(x)=5$. Assuming that $z^2=Q(u)$ is another hyperelliptic equation defining $C$, there exists
\[\left(
    \begin{array}{cc}
      a & b \\
      c & d \\
    \end{array}
  \right)\in\GL_2(K), e\in K^{\times}\] such that $\displaystyle x=\frac{a u+b}{c u+d}$, $\displaystyle y=\frac{ez}{(cu+d)^3}$.

  Assuming that $P(x)=a_0x^6+a_1x^5+a_2x^4+a_3x^3+a_4x^2+a_5x+a_6\in K[x]$, one may define the Igusa invariants (projective invariants), $J_{2i},1\le i\le 5$, associated to $P(x)$, see \cite{Liumodelesminimaux} for an explicit description of these invariants. One knows that $J_{2i}\in\Z[a_0,\ldots,a_6][1/2]$ is a homogeneous polynomial of degree $2i$ in the $a_i$'s. Moreover, one may define the invariants $I_4 = J_2^2 - 24J_4$ and $I_{12} = 2^{-2}(J_2^2J_4^2- 32J_4^3 - J_2^3J_6 + 36J_2J_4J_6 - 108J_6^2)$ which are homogeneous polynomials in the $a_i$'s of degree $4$ and $12$, respectively.

  If $J_{2i}',1\le i\le 5$, are the Igusa invariants associated to another hyperelliptic equation describing $C$, then there is an $a\in K^{\times}$ such that $J_{2i}'=a^{2i}J_{2i}$. Furthermore, if $C$ and $C'$ are two genus $2$ curves with Igusa invariants $J_{2i}$ and $J_{2i}'$ satisfying the latter equality, then $C$ and $C'$ are isomorphic over the algebraic closure of $K$.

  In \cite[\S2]{Liumodelesminimaux}, Liu introduced the following invariants (affine invariants) attached to $P(x)$:
  \begin{eqnarray*}
A_2 &=& -5a_1^2 + 12a_0a_2 \\
A_3 &=& 5a_1^3 + 9a_0(-2a_2a_1 + 3a_0a_3) \\
A_4 &=& -5a_1^4+ 24a_0(a_2a_1^2 - 3a_3a_0a_1 + 6a_4a_0^2)\\
A_5 &=& a_1^5+ 3a_0(-2a_2a_1^3 + 9a_0a_3a_1^2 - 36a_0^2a_4a_1 + 108a_0^3a_5)\\
B_2 &=& 2a_2^2 - 5a_1a_3 + 10a_0a_4.
\end{eqnarray*}
 One observes that $A_i,B_2\in\Z[a_0,\ldots,a_6]$.
 The invariants $A_i$ and $B_2$ are homogeneous polynomials in the $a_i$'s of degree $i$ and $2$, respectively.

\section{Models of genus two curves}
Throughout this article $R$ is a complete discrete valuation ring, $K$ its field of fractions, $\mathfrak{m}$ its maximal ideal generated by $t$, and $k$ its residue field. We assume moreover that $k$ is algebraically closed and $\Char k\ne2$. Let $\nu:K\to\Z\cup\{\infty\}$ be the normalized valuation of $K$, i.e., $\nu(t)=1$. The map $\alpha\mapsto\overline{\alpha}$ is the canonical homomorphism from $R$ to $k$.

For a scheme $X$ over $R$, $X_K$ and $X_k$ will denote the generic fiber and the special fiber of $X$, respectively.

\begin{Definition}
Let $C$ be a smooth projective curve over $K$, $R'$ a discrete valuation ring dominating $R$. We say that $\CC$ is a {\em stable model} for $C$ over $R'$ if $\CC$ is a stable curve over $R'$ whose generic fiber is isomorphic to $C\times_K K'$, where $K'$ is the fraction field of $R'$.
\end{Definition}
Therefore, one knows that the singular points of a stable model $\CC$ of $C$ over $R'$ are ordinary double points. In particular, if $P$ is a singular point of $\CC_k$, then the $\mathfrak{m}_P$-adic completion of $\mathcal{O}_{\CC,P}$ satisfies
\[ \widehat{\mathcal{O}}_{\CC,P}\simeq \widehat{R}[[u,v]]/(uv-\pi),\,\pi\in\mathfrak{m}_{K'}\setminus\{0\}.\] The {\em degree of singularity} of $P$ in $\CC$ is the integer $\nu_{K'}(\pi)$, where $\nu_{K'}$ is the normalized valuation of $K'$.

The following proposition, \cite[Proposition 4]{LiuCourbesstablesdegenre2}, shows that if the genus of $C$ is positive, then $C$ admits a stable model over a Galois extension of $K$.

\begin{Proposition}
\label{prop:fieldextension}
Suppose that $C$ is a smooth projective geometrically connected curve over $K$ of genus $g\ge 1$. There exists a Galois extension $L$ of $K$ such that for every finite extension $F$ of $K$, $C\times_K F$ admits a stable model if and only if $L\subseteq F$.
\end{Proposition}
\begin{Definition}
Let $C$ be a smooth projective geometrically connected curve of genus $g\ge 1$ over $K$. Let $L/K$ be the smallest Galois field extension over which $C$ admits a stable model $\CC$, see Proposition \ref{prop:fieldextension}. Then we say that $\CC$ is the {\em stable model} of $C$. Furthermore, the stable model $\CC$ of $C$ is unique up to isomorphism.
\end{Definition}

The special fiber $\CC_k$ of the stable model $\CC$ of $C$ is either smooth; irreducible with one or two double points; the union of two rational curves intersecting transversally in three points; or the union of two irreducible components intersecting in one and only one point. In the latter case, the two irreducible components are either smooth; singular; or one is smooth and the other component is singular. The following theorem, \cite[Th\'{e}or\`{e}me 1]{LiuCourbesstablesdegenre2}, gives explicit criteria for each of these possibilities of $\CC_k$ in terms of the invariants of $C$.

We set $\epsilon=1$ if $\Char k\ne 2,3$; $\epsilon =3$ if $\Char k=3$; and $\epsilon =4$ if $\Char k=2$. We define $I_2:=12^{-1}J_2,\,I_6:=J_6,\,I_8:=J_8$.
\begin{Theorem}
\label{Thm:singularity}
Let $C$ be a hyperelliptic curve defined over $K$ of genus $2$. Then
\begin{itemize}
\item[(I)] $\CC_k$ is smooth if and only if $J_{2i}^5J_{10}^{-i}\in R$ for every $i\le 5$;
\item[(II)] $\CC_k$ is irreducible with one double point if and only if $J_{2i}^6I_{12}^{-i}\in R$ for every $i\le 5$ and $J_{10}^6I_{12}^{-5}\in\mathfrak{m}$. The normalisation of $\CC_k$ is an elliptic curve with $j$-invariant $j=\overline{(I_4^3I_{12}^{-1})}$;
    \item[(III)] $C_k$ is irreducible with two double points if and only if $J_{2i}^2I_4^{-i}\in R$ for $i\le 5$, $J_{10}^2I_4^{-5}\in\mathfrak{m}$, $I_{12}I_4^{-3}\in\mathfrak{m}$, and either $J_4I_4^{-1}$ or $J_6^2I_4^{-3}$ is invertible in $R$;
        \item[(IV)] $\CC_k$ consists of two rational curves intersecting transversally in three points if and only if $J_{2i}^2I_4^{-i}\in\mathfrak{m}$ for $2\le i\le 5$;
        \item[(V*)] $\CC_k$ is the union of two irreducible components intersecting in one point if and only if \begin{eqnarray}
        I_4^{\epsilon}I_{2\epsilon}^{-2}\in\mathfrak{m},\; J_{10}^{\epsilon}I_{2\epsilon}^{-5}\in\mathfrak{m},\;  I_{12}^{\epsilon}I_{2\epsilon}^{-6}\in\mathfrak{m}
         \end{eqnarray} (which implies that $J_{2i}^{\epsilon}I_{2\epsilon}^{-i}\in R$ for $i\le5$). Furthermore,
  \item[(V)] Both components of $\CC_k$ are smooth if and only if: in addition to (1), $I_4^{3\epsilon}J_{10}^{-\epsilon}I_{2\epsilon}^{-1}\in R$, $I_{12}^{\epsilon}J_{10}^{-\epsilon}I_{2\epsilon}^{-1}\in R$. If $j_1$ and $j_2$ are the modular invariants of the components of $\CC_k$, then
      \[(j_1j_2)^{\epsilon}=\overline{(I_4^{3\epsilon}J_{10}^{-\epsilon}I_{2\epsilon}^{-1})},\;(j_1+j_2)^{\epsilon}=2^6.3^3+\overline{(I_{12}^{\epsilon}J_{10}^{-\epsilon}I_{2\epsilon}^{-1})};\]
  \item[(VI)] Only one of the two components of $\CC_k$ is smooth if and only if: in addition to (1), $I_4^3I_{12}^{-1}\in R, J_{10}^{\epsilon}I_{2\epsilon}I_{12}^{-\epsilon}\in\mathfrak{m}$. The modular invariant of the smooth component of $\CC_k$ is $j=\overline{(I_4^3I_{12}^{-1})}$;
  \item[(VII)] Both components of $\CC_k$ are singular if and only if: in addition to (1), $I_{12}I_4^{-3}\in\mathfrak{m}$, and $J_{10}^{\epsilon}I_{2\epsilon}I_4^{-3\epsilon}\in\mathfrak{m}$.
\end{itemize}
\end{Theorem}

The following proposition provides us with the degrees of singularity of the singular points of $\CC_k$.

\begin{Proposition}\cite[Proposition 2]{LiuCourbesstablesdegenre2}
\label{prop:degreeofsingularity}
Let $C$ be a smooth projective geometrically connected curve of genus $2$ over $K$. Let $L/K$ be the smallest Galois field extension over which $C$ has its stable model $\CC$, and $\nu_L$ its normalized valuation.
The following statements hold.
\begin{itemize}
\item[(I)] If $\CC_k$ is smooth, then the minimal regular model of $C$ is $\CC$.
\item[(II)] If $\CC_k$ is irreducible with one double point, the degree of singularity of its singular point is $e=\nu_L(J_{10}^6I_{12}^{-5})/6$.
\item[(III)] If $\CC_k$ is irreducible with two double points, the degrees of singularity, $e_1\le e_2$, of the singular points are given by
\[e_1=\min\{\nu_L(I_{12}I_4^{-3}),\nu_L(J_{10}^2I_4^{-5})\},\;e_2=\frac{1}{2}\nu_L(J_{10}^2I_4^{-5})-e_1.\]
\item[(IV)] If $\CC_k$ consists of two rational curves intersecting in three points, we assume $e_1\le e_2\le e_3$ are the degrees of singularity of the singular points of $\CC_k$. Set $l=\nu_L(J_{10}J_2^{-5}),n=\nu_L(I_{12}J_2^{-6})$ and $m=\nu_L(J_4J_2^{-2}).$ Then
    \[e_1=\min\{l/3,n/2,m\},\;e_2=\min\{(l-e_1)/2,n-e_1\},\;e_3=l-e_1-e_2.\]
\item[(V)] If $\CC_k$ is the union of two smooth irreducible components intersecting in one point, then the degree of singularity of the singular point of $\CC_k$ is $e=\nu_L(J_{10}^{\epsilon}I_{2\epsilon}^{-5})/12\epsilon$.
\item[(VI)] If $\CC_k$ is the union of one smooth and one singular irreducible component intersecting in one point, we set $e_0$ to be the degree of singularity of the point of intersection of the components of $\CC_k$, and $e_1$ the degree of singularity of the other singular point. Then \[e_0=\nu_L(I_{12}^{\epsilon}I_{2\epsilon}^{-6})/12\epsilon,\;e_1=\nu_L(J_{10}^{\epsilon}I_{2\epsilon}I_{12}^{-\epsilon})/\epsilon.\]
\item[(VII)] If $\CC_k$ is the union of two singular irreducible components intersecting in one point, we set $e_0$ to be the degree of singularity of the point of intersection of the two components of $\CC_k$ and $e_1\le e_2$ the degrees of singularity of the other singular points. Then $e_0=\nu_L(I_4^{\epsilon}I_{2\epsilon}^{-2})$,
    \[e_1=\min\{\nu_L(I_{12}I_4^{-3}),\nu_L(J_{10}^{\epsilon}I_{2\epsilon}I_4^{-3\epsilon})/2\epsilon\},\;e_2=\frac{1}{\epsilon}\nu_L(J_{10}^{\epsilon}I_{2\epsilon}I_4^{-3\epsilon})-e_1.\]
\end{itemize}
\end{Proposition}

\section{Stable models of quadratic twists}
Let $C$ be a smooth projective geometrically connected curve of genus $2$ defined over $K$. Let $y^2=P(x)$ be a defining polynomial of $C$. Let $K(\sqrt{D})$ be a quadratic extension of $K$ with associated character $\chi$. Let $C^{\chi}$ be the quadratic twist of $C$ by the character $\chi$. One knows that $C^{\chi}$ is defined by the hyperelliptic equation $y^2=DP(x)$. One may and will assume without loss of generality that $\nu(D)=0$ or $1$.

\begin{Proposition}
Let $C$ be a smooth projective geometrically connected curve of genus $2$ defined over $K$. Let $K(\sqrt{D})$ be a quadratic extension of $K$ with associated character $\chi$. If $\nu(D)=0$, then the minimal regular model $ X^{\chi}$ of $C^{\chi}$ is isomorphic to the minimal regular model $ X$ of $C$ over $R$.
\end{Proposition}
\begin{Proof}
  Since $\nu(D)=0$, one knows that $K(\sqrt D)$ is an unramified extension of $K$. Now, the proof follows from \cite[Proposition 10.1.17]{Liubook}.
\end{Proof}

The proposition above allows us to assume from now on that $\nu(D)=1$.
\begin{Lemma}
\label{Lem:IgusaofTwists}
   Let $J_{2i},1\le i\le 5$, $I_4$, $I_{12}$, $A_i,2\le i\le 5$, $B_2$ be the invariants attached to a hyperelliptic equation $y^2=P(x)$ defining a smooth projective geometrically connected genus $2$ curve over $K$, and $J_{2i}',1\le i\le 5$, $I_4'$, $I_{12}'$, $A_i',2\le i\le 5$, $B_2'$ the invariants attached to the hyperelliptic equation $y^2=DP(x)$. 
   Then one has $J_{2i}'=D^{2i}J_{2i}$, $I_4'=D^4I_4$, $I_{12}'=D^{12}I_{12}$, $A_{i}'=D^iA_i$, and $B_2'=D^2B_2$.
\end{Lemma}
\begin{Proof}
 One may assume that $P(x)=a_0x^6+a_1x^5+\ldots+a_6\in R[x]$. Now, $y^2=DP(x)=a_0'x^6+a_1'x^5+\ldots+a_6'$ where $a_i'=Da_i$. The result holds using the fact that the invariants $J_{i}$, $I_i$, $A_i$, and $B_i$ are homogeneous of degree $i$ in the $a_j$'s, see \S 2.
\end{Proof}

\begin{Lemma}
\label{Lem:StableModelofTwist}
Let $C$ and $C^{\chi}$ be as above. Let $\CC$ and $\CC^{\chi}$ be the stable models of $C$ and $C^{\chi}$ respectively. Then the following statements are true.
\begin{itemize}
\item[(I)] $\CC_k$ is smooth if and only if $\CC^{\chi}_k$ is smooth;
\item[(II)] $\CC_k$ is irreducible with one double point if and only if $\CC^{\chi}_k$ is irreducible with one double point;
    \item[(III)] $C_k$ is irreducible with two double points if and only if  $C^{\chi}_k$ is irreducible with two double points;
        \item[(IV)] $\CC_k$ consists of two rational curves intersecting transversally in three points if and only if $\CC^{\chi}_k$ consists of two rational curves intersecting transversally in three points;
        \item[(V*)] $\CC_k$ is the union of two irreducible components intersecting in one point if and only if $\CC^{\chi}_k$ is the union of two irreducible components intersecting in one point. Furthermore,
  \item[(V)] Both components of $\CC_k$ are smooth if and only if both components of $\CC^{\chi}_k$ are smooth;
  \item[(VI)] Only one of the two components of $\CC_k$ is smooth if and only if only only one of the two components of $\CC^{\chi}_k$ is smooth;
  \item[(VII)] Both components of $\CC_k$ are singular if and only if both components of $\CC^{\chi}_k$ are singular.
\end{itemize}
\end{Lemma}
\begin{Proof}
Let $C$ be defined by the hyperelliptic equation $y^2=P(x)$ and $C^{\chi}$ defined by $y^2=DP(x)$. If $J_{2i},\;1\le i\le 5,$ is an Igusa invariant attached to the hyperelliptic equation describing $C$, then $J_{2i}'$ will be an Igusa invariant of the hyperelliptic equation describing $C^{\chi}$. Similarly, $I_{2i}$ and $I_{2i}'$ are invariants of $C$ and $C^{\chi}$, respectively.

According to Theorem \ref{Thm:singularity}, one only needs to study a quotient of products of Igusa invariants powers. We prove (I) and (II) and the other cases are similar.
(I) One has that $\CC_k$ is smooth if and only if $J_{2i}^5J_{10}^{-i}\in R$ for every $i\le 5$. Since $J_{2i}$ is homogeneous of degree $2i$, one has $J_{2i}'^5J_{10}'^{-i}=J_{2i}^5J_{10}^{-i}\in R$. (II) $\CC_k$ is irreducible with one double point if and only if $J_{2i}^6I_{12}^{-i}\in R$ for every $i\le 5$ and $J_{10}^6I_{12}^{-5}\in\mathfrak{m}$, see Theorem \ref{Thm:singularity} (II). Since $J_{2i}^6I_{12}^{-i}$ and $J_{10}^6I_{12}^{-5}$ are quotients of invariants of the same degree, it follows that $J_{2i}'^6I_{12}'^{-i}=J_{2i}^6I_{12}^{-i}\in R$ for every $i\le 5$, and $J_{10}'^6I_{12}'^{-5}=J_{10}^6I_{12}^{-5}\in\mathfrak{m}$.
\end{Proof}

\section{Minimal regular models of quadratic twists}

Let $C$ be a smooth projective geometrically connected curve of genus $2$ defined over $K$. We let $\sigma$ be the hyperelliptic involution of $C$. It extends to an involution of the stable model $\CC$, which we will denote by $\sigma$ again. We set $\mathcal Z=\CC/\langle\sigma\rangle$ and $L$ the field extension of $K$ over which $C$ attains its stable model, see Proposition \ref{prop:fieldextension}. Then $\mathcal Z$ is a semi-stable curve over $R_L$, the ring of integers of $L$, whose generic fibre is isomorphic to $\PP_L^1$ and its special fibre is $\mathcal{Z}_k=\CC_k/\langle\sigma\rangle$. Let $f:\CC\to\mathcal Z$ be the canonical morphism, $\omega\in\mathcal{Z}_L$ the point corresponding to $x=\infty$, $\overline{\omega}\in\mathcal{Z}_k$ its specialization. We say $\overline{\omega}$ is ramified if $f$ is ramified above $\overline{\omega}$.

In \cite{NamikawaUeno}, the possible types for the special fiber of the minimal regular model of $C$ were introduced. In \cite{Liumodelesminimaux}, Liu displayed the possible reduction types when $\mathcal{Z}_k$ is either smooth, irreducible and singular, or not irreducible; and when the field extension $L/K$ is tamely ramified.
Moreover, he presented the possible reduction types for the minimal regular model of $C$ that appear when $L/K$ is a wildly ramified extension.

In this work, we follow the notation of Liu. In particular, the reduction type of the minimal regular model of $C$ will be given the same symbol as in \cite{Liumodelesminimaux}. Moreover, according to Lemma \ref{Lem:StableModelofTwist}, one knows that the special fibers of the stable models of the hyperelliptic curve $C$ and its quadratic twist have the same number of irreducible components and singular points. This reduces the number of possibilities that one needs to investigate in order to find the reduction type of the minimal regular model of the quadratic twist.

\begin{Remark}
\label{Rem1}
If $C^{\chi}$ is the quadratic twist of the curve $C$ defined over $K$ by $\chi$, then the quadratic twist of $C^{\chi}$ by $\chi$ is $C$ again. In particular, if $T$ is the reduction type of the minimal regular model of $C$ and $T'$ is the reduction type of the minimal regular model of $C^{\chi}$, then if $T'$ is the reduction type of the minimal regular model of some genus two curve, $T$ will be the reduction type of the minimal regular model of its quadratic twist by $\chi$. For example, in Theorem \ref{thm1} we will see that if $[\rm{IX}-3]$ ($[\rm{VIII}-1]$ respectively) is the reduction type of the minimal regular model of $C$, then the reduction type of the minimal regular model of the quadratic twist by $\chi$ is $[\rm{VIII}-1]$ ($[\rm{IX}-3]$ respectively).
\end{Remark}

\subsection{$\CC_k$ is smooth and $L/K$ is tamely ramified}
The following result can be found as \cite[Proposition 4.1.2 and Th\'{e}or\`{e}me 1]{Liumodelesminimaux}. It describes the ramification of the field extension $L/K$, see Proposition \ref{prop:fieldextension}, and its degree if the special fiber of the stable model is smooth.

We set $u_1=\nu(a_0^5J_{10}^{-1})$, $u_2=\nu(a_0^{10}J_{10}^{-1})$; and $v_1=\nu(A_5^{-2}J_{10})$, $v_2=\nu(A_5^{-6}J_{10}^{5})$, whereas $u_1',u_2',v_1',v_2'$ are the corresponding values for $C^{\chi}$.
\begin{Proposition}
\label{prop:smooth}
Assume that $\CC_k$ is smooth. The point $\overline{\omega}$ is ramified if and only if $A_5\ne 0$ and $a_0^{20}J_{10}A_5^{-6}\in\mathfrak{m}$. Moreover, the morphism $C\to \PP_K^1$ is ramified above some rational point $x_0\in \PP_K^1$.

Furthermore, the field extension $L/K$ is tamely ramified in each of the following situations:
\begin{itemize}
\item[(a)] $\Char k\ne 3,5$;
\item[(b)] $f:C\to \PP_K^1$ is ramified above two rational points in $\PP_K^1$;
\item[(c)] $\Char k =3$, $\overline{\omega}$ is ramified or $\CC_k\not\in \Gamma$ where $\Gamma$ is the set of isomorphis classes of the smooth proper curves over $k$ defined by $z^2=v^6+v^4+v^2+a,\;a\in k^{\times}$.
    \item[(d)] $\Char k= 5$, $\overline{\omega}$ is non-ramified or $J_{2i}^5J_{10}^{-i}\not\in\mathfrak{m}$ for some $i\le 3$.
\end{itemize}
In the case that $L/K$ is tamely ramified, we define $n,r,q$ as follows:
\begin{itemize}
\item[(a)] If $\overline{\omega}$ is non-ramified, $n$ is the least common denominator of $u_2/30$ and $u_1/10$, $r=nu_2/30$ and $q=n u_1/10$;
\item[(b)] If $\overline{\omega}$ is ramified, $n$ is the least common denominator of $v_1/20$ and $v_2/40$, $r=nv_1/20$ and $q=nv_2/40$. 
\end{itemize}
Then $[L:K]=n$.
\end{Proposition}

The integers $r$ and $q$ in the proposition above are used in the description of the action of $\Gal(L/K)$ on $\CC_k$, and will be used to determine the reduction type of the minimal regular model of the quadratic twist of $C$.

The field extension $L'$ is the field over which $C^{\chi}$ attains its stable model, and the degree of the extension is $n'$.

\begin{Lemma}
\label{lem:n-r-qsmooth}
Let $K(\sqrt{D})$, where $\nu(D)=1$, be a quadratic extension of $K$ with associated character $\chi$. Assume that $C$ is a smooth projective curve of genus two defined by $y^2=P(x)$ and $C^{\chi}$ is the quadratic twist of $C$ by $\chi$ defined by $y^2=DP(x)$. Let $n,q,r$ be the integers attached to $C$ defined in Proposition \ref{prop:smooth}, and $n',q',r'$ be the ones attached to $C^{\chi}$. We assume moreover that the special fiber of the stable model of $C$ is smooth. If $L/K$ is tamely ramified, then the following statements hold:
\begin{itemize}
\item[(a)] If $\overline{\omega}$ is non-ramified, then $\overline{\omega}'$ is non-ramified. Furthermore, $L'$ is a tamely ramified extension of $K$. The integer $n'$ is the least common denominator of $u_2/30$ and $u_1/10-1/2$, $nr'=n'r$ and $q'=n'(q/n-1/2)$;
\item[(b)] If $\overline{\omega}$ is ramified, then $\overline{\omega}'$ is ramified. Furthermore, $L'$ is a tamely ramified extension of $K$. The integer $n'$ is the least common denominator of $v_1/20$ and $v_2/40+1/2$, $nr'=n'r$ and $q'=n'(q/n+1/2)$.
\end{itemize}
\end{Lemma}
\begin{Proof}
In view of Proposition \ref{prop:smooth}, $\overline{\omega}'$ is ramified if and only if $A'_5\ne 0$ and $a_0^{20}J_{10}A_5^{-6}\in\mathfrak{m}$. The latter statement is equivalent to $\overline{\omega}$ being ramified since $A_5=D^{-5}A_5'\ne 0$ and $a_0^{20}J_{10}A_5^{-6}=a_0'^{20}J_{10}'A_5'^{-6}\in\mathfrak{m}$.

Again, we use the fact that the invariants of $C$ and $C^{\chi}$ are homogeneous in the $a_i$'s to evaluate $n',r'$ and $q'$. For (a), one has $n'$ is the least common denominator of $u_2'/30$ and $u_1'/10$, see Proposition \ref{prop:smooth}. Now $a_0'^{10}J_{10}'^{-1}=a_0^{10}J_{10}^{-1}$ and $a_0'^5J_{10}'^{-1}=a_0^5J_{10}^{-1}D^{-5}$. Now, the value for $n'$ follows from the fact that $\nu(ab)=\nu(a)+\nu(b)$ for $a,b\in K$ and $\nu(D)=1$. One has that $r'=n'u_2'/30=n'u_2/30=n'r/n$, and $q'=n'u_1'/10=n'\left(u_1-5\right)/10=n'(q/n-1/2)$.

The argument is similar for (b).
\end{Proof}

The tables in Theorem \ref{thm1} and Theorem \ref{thm:oneirreduciblecomponent} contain the reduction type of the minimal regular model of $C$, the positive integer $n$ which represents the degree of the field extension $L/K$ over which $C$ attains its stable model, the congruence classes of the two positive integers $r$ and $q$ mod $n$, see Proposition \ref{prop:smooth}, and the reduction type of the minimal regular model of the quadratic twist of $C$, together with the corresponding values $n'$, $r'$ mod $n'$ and $q'$ mod $n'$. In fact, determining the integers $n$, $r$ mod $n$ and $q$ mod $n$ yields the reduction type of the minimal regular model of $C$, see \cite[\S 8, Table 1, Table 2.1, Table 2.2, Table 2.3]{Liumodelesminimaux}. If the values of $q$ mod $n$ or $r$ mod $n$ do not appear in the table, then this means that there is no condition on these values when the corresponding reduction type occurs.

\begin{Theorem}
\label{thm1}
Let $C$ be a hyperelliptic curve defined over $K$. Assume that $L/K$ is tamely ramified and $\CC_k$ is smooth. Let $K(\sqrt{D})/K$ be a quadratic extension whose associated character is $\chi$, and $\nu(D)=1$. If $X$ and $X^{\chi}$ are the minimal regular models of the curves $C$ and its quadratic twist by $\chi$, $C^{\chi}$, then the reduction type of $X^{\chi}$ is given in the following table.

\begin{center}
\begin{tabular}{|c|c|c|c||c|c|c|c|}
  \hline
  type($X_k$) & $n$&$r$ mod $n$ & $q$ mod $n$ &type($X^{\chi}_k$)&$n'$& $r'$ mod $n'$&$q'$ mod $n'$\\
  \hline
   $[{\rm I}_{0-0-0}]$ &1&&& $[{\rm I}^*_{0-0-0}]$&2&0&\\
   \hline
   $[\rm{II}]$ &2& 1&&$[{\rm II}]$&2&1&\\
   \hline
   $[\rm{III}]$ &3& && $[\rm{IV}]$& 6 & 2 or 4&\\
   \hline
   $[\rm{VI}]$ &4&&& $[\rm{VI}]$&4&& \\
  \hline
  $[\rm{IX}-3]$ &5&1&& $[\rm{VIII}-1]$& 10&2&\\
  \hline
   $[\rm{IX}-1]$ &5&2&&$[\rm{VIII}-3]$&10&4& \\
   \hline
    $[\rm{IX}-4]$ &5&3&& $[\rm{VIII}-2]$&10&6&\\
    \hline
     $[\rm{IX}-2]$ &5&4&& $[\rm{VIII}-4]$&10&8&\\
     \hline
      $[\rm{V}]$ &6&1&0& $[\rm{V}^*]$& 6& 1& 3\\
      &6&5&3&&6&5&0\\
      \hline
         $[\rm{VII}^*]$ &8&&1 or 3&$[\rm{VII}]$& 8& & 5 or 7 \\
         \hline
       \end{tabular}
\end{center}
\end{Theorem}
\begin{Proof}
We assume that $C$ is given by the hyperelliptic equation $y^2=P(x)$, $\deg P(x)=5$ or $6$ and $C^{\chi}$ is given by $y^2=DP(x)$. Unless otherwise stated we will assume throughout the proof that  $\overline{\omega}$ is non-ramified, since the proof for $\overline{\omega}$ being ramified will be similar. We recall that if $\CC_k$ is smooth, then the special fiber $\CC^{\chi}_k$ of the stable model of $C^{\chi}$ is smooth, see Lemma \ref{Lem:StableModelofTwist}.

  We assume that the reduction type of $X_k$ is $[{\rm I}_{0-0-0}]$. Since $n=1$, $30\mid u_1$ and $10\mid u_2$, see Proposition \ref{prop:smooth}. Since $u'_1/30=u_1/30$ and $u'_2/10=u_2/10-1/2$, see Lemma \ref{lem:n-r-qsmooth}, it follows that $n'=2$. Moreover, $r'=n'r/n=2r\in2\Z$, i.e., $r'=0$ mod $2$. One has that the type of $X^{\chi}_k$ is $[{\rm I}^*_{0-0-0}]$.

If $X_k$ has reduction type $[{\rm II}]$, then $15\mid u_1$ and $2\nmid u_1$ since $r=2u_1/30=1$ mod $2$. Therefore, $n'=2$ as $\displaystyle u_1/30\in\frac{1}{2}\Z$ and the least denominator of $u_2/10-1/2$ is either $1$ or $2$. Furthermore, $r'=n'r/n=r=1$ mod $2$. One obtains that the reduction type of $X^{\chi}_k$ is $[{\rm II}]$.

If $X$ has reduction type $[{\rm III}]$, then $10\mid u_1,\;3\nmid u_1$ and $10\mid u_2$ since $n=3$. Now, $\displaystyle \frac{u_1}{30}\in\frac{1}{3}\Z$ and $\displaystyle \frac{u_2}{10}-\frac{1}{2}\in\frac{1}{2}\Z$ which implies that $n'=6$, see Lemma \ref{lem:n-r-qsmooth}. One moreover has $\displaystyle r'=6r/3=2r\in 2\Z$, it follows that $r'= 2$ or $4$ mod $6$. Thus, the type of $X^{\chi}_k$ is $[{\rm IV}]$. One observes that if $\overline{\omega}$ is ramified, then the reduction type of $X$ cannot be $[{\rm III}]$ since $3$ does not divide the denominator of $v_1/20$ and $v_2/40$.

If $X$ has reduction type $[{\rm VI}]$, then $n=4$. If $\overline{\omega}$ is non-ramified, then the least common denominator of $u_1/30$ and $u_2/10$ cannot be $4$. Therefore, $\overline{\omega}$ is ramified. One has $n'$ is the least common denominator of $\displaystyle \nu(v_1)/20$ and $\nu(v_2)/40+1/2$, see Lemma \ref{lem:n-r-qsmooth}. If the least denominator of $\nu(v_1)/20$ is $4$, then $n'=4$; otherwise the least denominator of $\nu(v_2)/40$ is $4$ and so $\nu(v_2)/10=1$ or $3$ mod $4$, thus $\nu(v_2)/10+2= 1$ or $3$ mod $4$, which implies that $n'=4$. Therefore, the type of $X^{\chi}_k$ is $[{\rm VI}]$.

Assuming that $n=5$, and $\overline{\omega}$ is non-ramified, one has $5$ is the least common denominator of $ u_1/30$ and $u_2/10$. Now, $n'$ is the least common denominator of $\displaystyle u_1/30$ and $u_2/10-1/2$. This yields that $n'=10$. Moreover, $r'=2r$. Therefore, if the reduction type of $X$ is $[\textrm{IX}-i]$, $i=3, 1, 4,2$, then $r=1,2,3,4$ mod $5$, respectively, and so $r'=2, 4, 6, 8$ mod $10$, respectively, i.e., the reduction type of $X^{\chi}$ is $[\textrm{VIII}-j]$, $j=1,3,2,4$, respectively.

Assuming that $n=6$, $\overline{\omega}$ must be non-ramified. One has $5\mid u_1$ and $5\mid u_2$. If $r:=6u_1/30=1$ or $5$ mod $6$, then $r=5u_1$ mod $6$, and $u_1=5$ or $1$ mod $6$, respectively. In particular, $2\nmid u_1$, hence $\displaystyle u_1/30\in\frac{1}{6}\Z$. In this case, the least common denominator, $n'$, of $u_1/30$ and $u_2/10-1/2$ is $6$, moreover, $r'=r$ and $\displaystyle q'= q-3$, see Lemma \ref{lem:n-r-qsmooth}. This implies that if $X$ has reduction type $[{\rm V}]$ or $[{\rm V}^*]$, then the reduction type of $X^{\chi}$ is $[{\rm V}^*]$ or $[{\rm V}]$, respectively.

If $n=8$, then $\overline{\omega}$ must be ramified. Now, the least denominator of $v_2/40$ is $8$, in particular, $v_2/5=1,3,5$ or $7$ mod $8$. It follows that $\displaystyle v_2/40+1/2\in\frac{1}{8}\Z$. Therefore, $n'=8$, $r'=r$, and $q'=q+4$, Lemma \ref{lem:n-r-qsmooth}. If $q=1$ or $3$ mod $8$ ($5$ or $7$ mod $8$), then $q'=5$ or $7$ mod $8$ ($1$ or $3$ mod $8$). More specifically, if the reduction type of $X$ is $[\rm{VII}^*]$ or $[\rm{VII}]$, then the reduction type of $X^{\chi}$ is $[\rm{VII}]$ or $[\rm{VII}^*]$, respectively. Hence this covers all the possible reduction types when $L/K$ is tamely ramified and $\CC_k$ is smooth.
\end{Proof}

\begin{Remark}
One observes that there are genus two curves whose twists have the same reduction type as the curve itself. For example, those genus two curves with reduction type ${\rm II}$ and ${\rm VI}$. This is a phenomenon that does not occur for elliptic curves.
\end{Remark}
\subsection{$\CC_k$ is singular, $\CC_k/\langle\sigma\rangle$ is irreducible, and $L/K$ is tamely ramified}

In this section we assume that $\CC_k/\langle\sigma\rangle$ is irreducible and $\CC_k$ is singular. In particular, $\CC_k$ is irreducible with one or two double point; or $\CC_k$ consists of two rational curves intersecting transversally in three ordinary double points.

Let $C_{000}$ be the stable curve over $k$ consisting of two rational curves intersecting in three points. We set:

\begin{align*}J_{12}=\left\{\begin{array}{ll}
I_{12}  & \textrm{if $\CC_k$ has a singular point} \\
I_4^3 & \textrm{if $\CC_k$ is irreducible and rational}\\
J_2^6 & \textrm{if $\CC_k=C_{000}$} \end{array}\right.\end{align*}

The following proposition, see \cite[Proposition 4.2.1, Proposition 4.2.2 and Th\'{e}or\`{e}me 2]{Liumodelesminimaux}, describes the ramification of $f:\CC\to \mathcal Z$ over $\overline{\omega}$, the cases when $L/K$ is tamely ramified, and the field extension $L/K$ over which $C$ attains its stable model.

We set $u_1=\nu(a_0^{12}J_{12}^{-1})$, $u_2=\nu(a_0^6J_{12}^{-1})$.

\begin{Proposition}
\label{prop:irreduciblesingular}
Assume that $\CC_k$ is singular and $\CC_k/\langle\sigma\rangle$ is irreducible.
 \begin{itemize}
 \item[(a)] The point $\overline{\omega}$ is non-ramified if and only if $a_0^{-6}B_2^9J_{12}^{-1},a_0^{-120}A_5^{36}J_{12}^{-5}\in R$;
 \item[(b)] The point $\overline{\omega}$ is ramified and $f^{-1}(\overline{\omega})$ is a regular point if and only if $a_0^{120}A_5^{-36}J_{12}^{5}\in \mathfrak{m},\;B_2^{60}A_5^{-12}J_{12}^{-5}\in R$;
     \item[(c)] The point $\overline{\omega}$ is ramified and $f^{-1}(\overline{\omega})$ is a singular point if and only if $a_0^6B_2^{-9}J_{12}\in\mathfrak{m}$, and $B_2^{-60}A_5^{12}J_{12}^5\in\mathfrak{m}$.
\end{itemize}
Furthermore, the field extension $L/K$ is tamely ramified in each of the following cases:
\begin{itemize}
\item[(a)] $\Char k\ne 3$ or $\CC_k\ne C_{000}$;
\item[(b)] $\Char k=3$, $\CC_k=C_{000}$ and $\overline{\omega}$ is ramified.
\end{itemize}
In the case that $L/K$ is tamely ramified, we define $n,r,q$ as follows:
\begin{itemize}
\item[(a)] If $\overline{\omega}$ is non-ramified, $n$ is the least common denominator of $u_1/36$ and $u_2/12$, $r=nu_1/36$ and $q=nu_2/12$;
\item[(b)] If $\overline{\omega}$ is ramified and $f^{-1}(\overline{\omega})$ is a regular point, $n$ is the least denominator of $\nu(A_5^{36}J_{12}^{-25})/240$, $q=n\nu(A_5^{36}J_{12}^{-25})/240$ and $r=-2q$;
\item[(c)] If $\overline{\omega}$ is ramified and $f^{-1}(\overline{\omega})$ is singular, $n$ is the least common denominator of $\nu(B_2^{-6}J_{12})/12$ and $\nu(B_2^{-9}J_{12})/12$, $r=n\nu(B_2^{-6}J_{12})/12$ and $q=n\nu(B_2^{-9}J_{12})/12$.
\end{itemize}
Then $[L:K]=n$.
\end{Proposition}

The following lemma introduces the degree of the field extension, $n'$, over which $C^{\chi}$ admits a stable model, and the integers $r'$, $q'$ attached to $C^{\chi}$.

\begin{Lemma}
\label{lem:n-r-qirrdeuciblesingular}
Assume that $\CC_k$ is singular and $\CC_k/\langle\sigma\rangle$ is irreducible. Assume that $C$ is defined by $y^2=P(x)$ and $C^{\chi}$ is defined by $y^2=DP(x)$ where $\nu(D)=1$. We assume moreover that $L/K$ is tamely ramified. The following statements hold:
\begin{itemize}
\item[(a)] If $\overline{\omega}$ is non-ramified, then $\overline{\omega}'$ is non-ramified. Furthermore, $L'$ is a tamely ramified extension of $K$. The integer $n'$ is the least common denominator of $u_1/36$ and $u_2/12-1/2$, $nr'=n'r$ and $q'=n'(q/n-1/2)$;
\item[(b)] If $\overline{\omega}$ is ramified and $f^{-1}(\overline{\omega})$ is a regular point, then $\overline{\omega}'$ is ramified and $f'^{-1}(\overline{\omega'})$ is regular. Furthermore, $L'$ is a tamely ramified extension of $K$. The integer $n'$ is the least denominator of $\nu(A_5^{36}J_{12}^{-25})/240-1/2$, $q'=n'(q/n-1/2)$ and $r'=-2q'$;
\item[(c)] If $\overline{\omega}$ is ramified and $f^{-1}(\overline{\omega})$ is singular, then $\overline{\omega}'$ is ramified and $f'^{-1}(\overline{\omega}')$ is singular. Furthermore, $L'$ is a tamely ramified extension of $K$. The integer $n'$ is the least common denominator of $\nu(B_2^{-6}J_{12})/12$ and $\nu(B_2^{-9}J_{12})/12-1/2$, $nr'=n'r$ and $q'=n'(q/n-1/2)$.
\end{itemize}
\end{Lemma}
\begin{Proof}
This follows from Proposition \ref{prop:irreduciblesingular} and the fact that $a_0,B_2, A_5,J_{12}$ are homogeneous polynomials in the $a_i$'s of degrees $1,2,5,12$, respectively.
\end{Proof}

\begin{Theorem}
\label{thm:oneirreduciblecomponent}
Let $C$ be a hyperelliptic curve defined over $K$. Let $K(\sqrt{D})/K$ be a quadratic extension whose associated character is $\chi$, and $\nu(D)=1$. Let $X$ and $X^{\chi}$ be the minimal regular models of the curves $C$ and its quadratic twist by $\chi$, $C^{\chi}$.
\begin{itemize}
\item[a)] If $L/K$ is tamely ramified and $\CC_k$ consists of one irreducible component with a unique double point, then the reduction type of $X^{\chi}$ is given in the following table.
\vskip10pt
\begin{center}
\begin{tabular}{|c|c|c|c||c|c|c|c|}
  \hline
  type($X_k$) &$n$& $r$ mod $n$ & $q$ mod $n$ &type($X^{\chi}_k$)& $n'$& $r'$ mod $n'$& $q'$ mod $n'$ \\
  \hline
   $[{\rm I}_{d-0-0}]$ & 1& & & $[{\rm I}^*_{d-0-0}]$& 2 &0 &\\
   \hline
  $[{\rm II}_{d/2-0}]$ &2&1&1& $[{\rm II}^*_{d/2-0}]$&2&1&0\\
  \hline
  $[{\rm IV}-{\rm II}_{(d-2)/3}]$ &3&1&& $[{\rm II}^*-{\rm II}^*_{(d-2)/3}]$&6&2&\\
  \hline
   $[{\rm IV}^*-{\rm II}_{(d-1)/3}]$ &3&2&&$[{\rm II}-{\rm II}^*_{(d-1)/3}]$ &6&4&\\
   \hline
    $[{\rm III}-{\rm II}_{(d-2)/4}]$ &4&1&1& $[{\rm III}^*-{\rm II}^*_{(d-2)/4}]$&4&1&3\\
    \hline
      $[{\rm III}-{\rm II}^*_{(d-2)/4}]$ &4&3&1 &$[{\rm III}^*-{\rm II}_{(d-2)/4}]$&4&3&3 \\
      \hline
\end{tabular}
\end{center}
\item[b)] If $L/K$ is tamely ramified and $\CC_k$ consists of one irreducible component with exactly two ordinary double points, then the reduction type of $X^{\chi}$ is given in the following table.
\vskip10pt
\begin{center}
\begin{tabular}{|c|c|c|c||c|}
  \hline
  type($X_k$)& $n$ &$r$ mod $n$  & $f^{-1}(\omega)$ &type($X^{\chi}_k$)\\
  \hline
   $[{\rm I}_{d_1-d_2-0}]$ & 1& && $[{\rm I}^*_{d_1-d_2-0}]$\\
   \hline
   $[{\rm I}^*_{d_1/2-d_2/2-0}]$&2 &0& &$[{\rm I}_{d_1/2-d_2/2-0}]$\\
   \hline
   $[2{\rm I}_{d_1}-0]$ &2&1&regular &$[2{\rm I}_{d_1}-0]$ \\
   \hline
   $[{\rm II}_{d_1/2-d_2/2}]$ &2&1 & singular&$[{\rm II}_{d_1/2-d_2/2}]$  \\
   \hline
   $[{\rm III}_{d_1/2}]$ & 4&&  &$[{\rm III}_{d_1/2}]$\\
   \hline
   \end{tabular}
   \end{center}
\item[c)] If $L/K$ is tamely ramified and $\CC_k$ is the union of two rational curves intersecting transversally in three ordinary double points, then the reduction type of $X^{\chi}$ is given in the following table.
\vskip10pt
\begin{center}
\begin{tabular}{|c|c|c|c||c|c|c|c|}
  \hline
  type($X_k$) &$n$& $r$ mod $n$ &$q$ mod $n$ &type($X^{\chi}_k$)&$n'$&$r'$ mod $n'$& $q'$ mod $n'$\\
  \hline
   $[{\rm I}_{d_1-d_2-d_3}]$ &1&&&$[{\rm I}^*_{d_1-d_2-d_3}]$&2&0& \\
   \hline
$[{\rm II}_{e_1/2-e_2}]$& 2&1&1&$[{\rm II}^*_{e_1/2-e_2}]$&2&1&0\\
\hline
$[{\rm III}_{d_1}]$& 3&&&$[{\rm III}^*_{d_1}]$&6&&\\
\hline
\end{tabular}
\end{center}
\end{itemize}
   \end{Theorem}
  \begin{Proof}
 a) According to Lemma \ref{lem:n-r-qirrdeuciblesingular}, we have three subcases to consider: $\overline{\omega}$ is non-ramified, $f^{-1}(\overline{\omega})$ is regular, or $\overline{\omega}$ is ramified and $f^{-1}(\overline{\omega})$ is singular. Unless otherwise stated, we will assume that $\overline{\omega}$ is non-ramified since the proofs for the other two subcases will be similar. We recall that the degree of singularity of the unique double point in $\CC_k$ is given by $\nu_L(J_{10}^6I_{12}^{-5})/6$, see Proposition \ref{prop:degreeofsingularity} (II). Moreover, $\CC^{\chi}_k$ consists of one irreducible component with a unique double point.

If $X$ has reduction type  $[{\rm I}_{d-0-0}]$, then this implies that the degree of singularity of the ordinary double point in $\CC_k$ is $d$. In view of Proposition \ref{prop:irreduciblesingular}, since $n=1$, one has $\displaystyle 36\mid u_1$ and $12 \mid u_2$. Now, $n'$ is the least common denominator of $u_1/36$ and $u_2/12-1/2$, see Lemma \ref{lem:n-r-qirrdeuciblesingular}. This yields $n'=2$, $r'=2r= 0$ mod $2$. Therefore, the reduction type of $X^{\chi}$ is $[{\rm I}^*_{d'/2-0-0}]$ where $d'$ is the degree of singularity of the double point of $\CC^{\chi}_k$. One observes that $d'=\nu_{L'}(J_{10}'^6I_{12}'^{-5})/6=2\nu(J_{10}^6I_{12}^{-5})/6=2d$.

When $n=2$ and $r=1$ mod $2$, one has $\displaystyle u_1/36\in\frac{1}{2}\Z$. Since $n'$ is the least common denominator of $u_1/36$ and $u_2/12-1/2$, it follows that $n'=2$, $r'=r$, and the degree of singularity of the double point in $\CC^{\chi}_k$ is $d'=d$. Moreover, one has $q'=q-1$. Thus, if $q=1$ mod $2$, then the reduction type of $X$ is $[{\rm II}_{d/2-0}]$ and $q'=0$ mod $2$, hence the reduction type of $X^{\chi}$ is $[{\rm II}^*_{d/2-0}]$. One observes that when $r=1$ mod $2$, if $\overline{\omega}$ is ramified, then $f^{-1}(\overline{\omega})$ cannot be regular as according to Proposition \ref{prop:irreduciblesingular}, $r$ must be even.

If $n=3$, then the least common denominator of $u_1/36$ and $u_2/12$ is $3$. Since $n'$ is the least common denominator of $u_1/36$ and $u_2/12-1/2$, it follows that $n'=6$ and $r'=2r$. If $r=1$ mod $3$ (the reduction type of $X$ is $[{\rm IV}-{\rm II}_{(d-2)/3}]$),  then $r'=2$ mod $6$ and the reduction type of $X^{\chi}$ is $[{\rm II}^*-{\rm II}^*_{(d'-4)/6}]$ where $d'=2d$. If $r=2$ mod $3$ (the reduction type of $X$ is $[{\rm IV}^*-{\rm II}_{(d-1)/3}]$), then $r'= 4$ mod $6$ and the reduction type of $X^{\chi}$ is $[{\rm II}-{\rm II}^*_{(d'-2)/6}]$ where $d'=2d$.

 If $n=4$, then the least common denominator of $u_1/36$ and $u_2/12-1/2$ is $n'=4$. Moreover, $r'=r$, $q'=q-2$, and the degree of singularity $d'=d$.

The proofs of b) and c) follow the same lines of the proof of a), therefore, they will be omitted.
\end{Proof}
  The integers $d_1$ and $d_2$ appearing in the reduction types of Theorem \ref{thm:oneirreduciblecomponent} b) are the degrees of singularity of the two ordinary double points in $\CC_k$. The degrees of singularity of the double points is given in Proposition \ref{prop:degreeofsingularity} (III).
The degrees of singularity of the two double points in $\CC_k$ when the reduction type of $X$ is either $[2{\rm I}_{d_1}-0]$ or $[{\rm III}_{d_1/2}]$ are both equal to $d_1$.

The integers $d_1,d_2$ and $d_3$ in the reduction types appearing in Theorem \ref{thm:oneirreduciblecomponent} c) are the degrees of singularity of the three ordinary double points, and can be evaluated using Proposition \ref{prop:degreeofsingularity} (IV).
When the reduction type of $X$ is either $[{\rm II}^*_{e_1/2-e_2}]$ or $[{\rm II}_{e_1/2-e_2}]$, then exactly two of the ordinary double points in $\CC_k$ have the same degree of singularity $e_1$, and the degree of singularity of the third ordinary double point is $e_2$. If the reduction type of $X$ is either $[{\rm III}_{d_1}]$ or $[{\rm III}^*_{d_1/2}]$, then the three ordinary double points of $\CC_k$ have the same degree of singularity $d_1$.

\subsection{$\CC_k/\langle\sigma\rangle$ is not irreducible and $L/K$ is tamely ramified}
Now, we assume that $\CC_k/\langle\sigma\rangle$ is not irreducible where the field $L$ over which $C$ admits a stable model is a tamely ramified extension.
Assuming that $\Char k\ne 3$ and $\mathcal{Z}_k$ is the union of two projective curves intersecting in one point. The possible divisors of $[L:K]$ are $2$ and $3$. It follows that $L/K$ is tamely ramified. Letting $E_1$ and $E_2$ be the irreducible components of $\CC_k$, with $\overline{\omega}\in f(E_1)$, we set:
\begin{align*}d_K=\left\{\begin{array}{ll}
\nu(J_{10}J_2^{-5})/12  & \textrm{if $E_1,E_2$ are smooth} \\
\nu(I_{12}J_2^{-6})/12 & \textrm{if $\CC_k$ has a unique smooth component}\\
\nu(I_4J_2^{-2})/4& \textrm{if $E_1,E_2$ are singular} \end{array}\right.\end{align*}
The degree of singularity of the point of intersection $E_1\cap E_2$ in $\CC$ is $d=[L:K]d_K$. It follows that the degree of singularity of $f(E_1\cap E_2)$ in $\mathcal{Z}$ is $2d$.

The following proposition, see \cite[Proposition 4.3.1 and Th\'{e}or\`{e}me 3]{Liumodelesminimaux}, summarizes the behavior of the stable model in that case.

\begin{Proposition}
\label{prop:notirreducible}
Assume $\Char k\ne 3$. Assume that $\CC_k/\langle\sigma\rangle$ is not irreducible.
 \begin{itemize}
 \item[(a)] The point $\overline{\omega}$ is non-ramified if and only if $a_0^{-2}B_2^3J_{2}^{-2}\in R,\;a_0^{-4}A_3^{2}J_{2}^{-1},a_0^{-20}A_5^6J_2^{-5}\in R$ and at least one of the two latter elements are invertible in $R$;
 \item[(b)] $f^{-1}(\overline{\omega})$ is a regular point if and only if $a_0^{20}A_5^{-6}J_{2}^{5}\in \mathfrak{m}$ and $B_2^{10}A_5^{-2}J_{2}^{-5}\in R$;
     \item[(c)] $\overline{\omega}$ is regular and $f^{-1}(\overline{\omega})$ is a singular point if and only if $a_0^2B_2^{-3}J_{2}^2\in\mathfrak{m}$ and $B_2^{-10}A_5^{2}J_{2}^5\in\mathfrak{m}$;
          \item[(d)] $\overline{\omega}$ is singular if and only if $a_0^{-2}B_2^{3}J_{2}^{-2}\in R$ and $a_0^{-4}A_3^2J_2^{-1},a_0^{-20}A_5^6J_2^{-5}\in\mathfrak{m}$.
\end{itemize}
Assume that $2\mid \nu(J_2)$, then one has
\begin{itemize}
\item[(a)] If $\overline{\omega}$ is non-ramified, $n$ is the smallest common denominator of $d_K$ and $\nu(a_0J_{2})/6$, $r=n\nu(a_0J_{2})/6$;
\item[(b)] If $\overline{\omega}$ is regular and ramified, and $f^{-1}(\overline{\omega})$ is a regular point, then $n$ is the least common denominator of $d_K$ and $\nu(A_5^{2}J_{2})/8$, $r=n\nu(A_5^2J_2)/8$;
\item[(c)] If $\overline{\omega}$ is regular such that $f^{-1}(\overline{\omega})$ is singular, then $n$ is the least common denominator of $d_K$ and $\nu(B_2)/4$, $r=n\nu(B_2)/4$;
    \item[(d)] If $\overline{\omega}$ is singular, $n$ is the least common denominator of $d_K$ and $r_K$, $r=nr_K$, where
    \[r_K=\nu(a_0)/2+\min\{d_K/2,\nu(A_2^{-3}A_3^2)/8,\nu(A_2^{-5}(A_2A_3-3A_5)^2)/12\}\in\Q.\]
\end{itemize}
Then $[L:K]=n$ and $d=nd_K$.

If $2\nmid \nu(J_2)$, then $d_K+\nu(a_0)=2r_K$ and $[L:K]=2m$, where $m$ is the least denominator of $d_K$, and $r=md_K$.
\end{Proposition}
\begin{Lemma}
\label{lem:n-r-qnotirreducible}
Assume $\Char k\ne 3$. Assume that $\CC_k/\langle\sigma\rangle$ is not irreducible. Then $L'/K$ is tamely ramified. Moreover,
 \begin{itemize}
 \item[(a)] If the point $\overline{\omega}$ is non-ramified, then $\overline{\omega}'$ is non-ramified;
 \item[(b)] If $f^{-1}(\overline{\omega})$ is a regular point, then $f'^{-1}(\overline{\omega}')$ is regular;
     \item[(c)] If $\overline{\omega}$ is regular and $f^{-1}(\overline{\omega})$ is a singular point, then $\overline{\omega}'$ is regular and $f'^{-1}(\overline{\omega}')$ is singular;
          \item[(d)] If $\overline{\omega}$ is singular, then $\overline{\omega}'$ is singular.
\end{itemize}
Moreover, if $2\mid \nu(J_2)$, then $2\mid \nu(J_2')$, and the following statements hold:
\begin{itemize}
\item[(a)] If $\overline{\omega}$ is non-ramified, $n'$ is the least common denominator of $d_K$ and $\nu(a_0J_{2})/6+1/2$, $nr'=n'(r+n/2)$;
\item[(b)] If $\overline{\omega}$ is regular and ramified, and $f^{-1}(\overline{\omega})$ is a regular point, then $n'$ is the least common denominator of $d_K$ and $\nu(A_5^{2}J_{2})/8+3/2$, $nr'=n'(r+3n/2)$;
\item[(c)] If $\overline{\omega}$ is regular such that $f^{-1}(\overline{\omega})$ is singular, then $n'$ is the least common denominator of $d_K$ and $\nu(B_2)/4+1/2$, $nr'=n'(r+n/2)$;
    \item[(d)] If $\overline{\omega}$ is singular, $n$ is the least common denominator of $d_K$ and $r_K+1/2$, $nr'=n'(r+n/2)$.
\end{itemize}
Moreover, $nd'=n'd$.

If $2\nmid \nu(J_2)$, then $2\nmid \nu(J_2')$. Furthermore, $[L':K]=[L:K]=2m$, where $m$ is the least denominator of $d_K$, and $r'=r$.
\end{Lemma}
\begin{Proof}
We assume that $C$ is defined by $y^2=P(x)=a_0x^6+\ldots+a_6\in R[x]$, and $C'$ is defined by $y^2=DP(x)$ and $\nu(D)=1$.

According to Proposition \ref{prop:notirreducible}, $\overline{\omega}'$ is non-ramified if $a_0'^{-2}B_2'^3J_{2}'^{-2}\in R$, $a_0'^{-4}A_3'^{2}J_{2}'^{-1},$ $a_0'^{-20}A_5'^6J_2'^{-5}\in R$ and at least one of the two latter elements are invertible in $R$. Observing that $a_0'^{-2}B_2'^3J_{2}'^{-2}=a_0^{-2}B_2^3J_{2}^{-2}$, $a_0'^{-4}A_3'^{2}J_{2}'^{-1}=a_0^{-4}A_3^{2}J_{2}^{-1}$ and $a_0'^{-20}A_5'^6J_2'^{-5}=a_0^{-20}A_5^6J_2^{-5}$, it follows that $a_0^{-2}B_2^3J_{2}^{-2},a_0^{-4}A_3^{2}J_{2}^{-1},a_0^{-20}A_5^6J_2^{-5}\in R$ and at least one of the two latter elements are invertible in $R$. This implies that $\overline{\omega}'$ is non-ramified if $\overline{\omega}$ is non-ramified. The proofs of (b), (c), and (d) are similar.

 Since $J_2'=D^2J_2$, one has $2\mid(\nu(D^2)+\nu(J_2))=2+\nu(J_2)$ if and only if $2\mid \nu(J_2)$.

If $\overline{\omega}$ is ramified, then $\overline{\omega}'$ is ramified and $n'$ is the least common denominator of $d'_K=d_K$ and $\nu(a_0'J_2')/6=\nu(a_0J_2)/6+1/2$. Moreover, $r'=n'\nu(a_0'J_2')/6=n'(\nu(a_0J_2)+3)/6$. The same argument holds for the other subcases.
\end{Proof}
The tables in the following theorem contain the reduction type of the minimal regular model of $C$, the positive integer $n$ which represents the degree of the field extension $L/K$ over which $C$ attains its stable model, the congruence classes of the two positive integers $d$ mod $n$ and $r$ mod $n$, see Proposition \ref{prop:notirreducible} and Lemma \ref{lem:n-r-qnotirreducible}, and the reduction type of the minimal regular model of the quadratic twist of $C$, together with the corresponding values $n'$, $d'$ mod $n'$ and $r'$ mod $n'$. In fact, determining the integers $n$, $d$ mod $n$ and $r$ mod $n$ yields the reduction type of the minimal regular model of $C$ when $\CC_k/\langle\sigma\rangle$ is not irreducible, see \cite[\S 8, Table 3.1, Table 3.2, Table 3.3, Table 3.4]{Liumodelesminimaux}. If the values of $d$ mod $n$ or $r$ mod $n$ do not appear in the table, then this means that there is no condition on these values when the corresponding reduction type occurs.

\begin{Theorem}
\label{thm:twoellipticcurves}
Let $C$ be a hyperelliptic curve defined over $K$. Let $K(\sqrt{D})/K$ be a quadratic extension whose associated character is $\chi$, and $\nu(D)=1$. Let $X$ and $X^{\chi}$ be the minimal regular models of the curves $C$ and its quadratic twist by $\chi$, $C^{\chi}$.
\begin{itemize}
\item[a)] Assume that $L/K$ is tamely ramified and $\CC_k$ is the union of two elliptic curves intersecting in one point.
\begin{itemize}
\item[i)] If $2\mid\nu(J_2)$, then the reduction type of $X^{\chi}$ is given in the following table.
\end{itemize}
\end{itemize}
{\footnotesize\begin{tabular}{|c|c|c|c||c|c|c|c|}
  \hline
  type($X_k$) &$n$ & $d$ mod $n$ & $r$ mod $n$ &type($X^{\chi}_k$)& $n'$& $d'$ mod $n'$ & $r'$ mod $n'$\\
  \hline
   $[{\rm I}_0-{\rm I}_0-d]$ & 1 & & & $[{\rm I}^*_0-{\rm I}^*_0-(d-1)]$&2&0&\\
   \hline
   $[{\rm I}_0-{\rm I}^*_0-(d-1)/2]$ &2& 1 & &$[{\rm I}_0-{\rm I}^*_0-(d-1)/2]$&2&1& \\
   \hline
   $[{\rm IV}-{\rm IV}^*-(d-3)/3]$& 3 & 0&  &$[{\rm II}-{\rm II}^*-(d-3)/3]$&6&0&\\
   \hline
   $[{\rm I}_0-{\rm IV}-(d-1)/3]$ & 3 & 1& 0 or 1&$[{\rm I}^*_0-{\rm II}^*-(d-4)/3]$&6&2&3 or 5\\
   \hline
   $[{\rm IV}^*-{\rm IV}^*-(d-4)/3]$ &3 & 1 & 2& $[{\rm II}-{\rm II}-(d-1)/3]$&6&2&1\\
   \hline
   $[{\rm I}_0-{\rm IV}^*-(d-2)/3]$ & 3 & 2& 0 or 2& $[{\rm I}_0^*-{\rm II}-(d-2)/3]$&6&4&3 or 1\\
   \hline
   $[{\rm IV}-{\rm IV}-(d-2)/3]$ &3& 2& 1&$[{\rm II}^*-{\rm II}^*-(d-5)/3]$&6&4&5\\
   \hline
   $[{\rm III}-{\rm III}^*-(d-4)/4]$ &4&0&&$[{\rm III}-{\rm III}^*-(d-4)/4]$&4&0&\\
   \hline
   $[{\rm I}_0-{\rm III}-(d-1)/4]$ &4&1& 0 or 1& $[{\rm I}^*_0-{\rm III}^*-(d-5)/4]$&4&1& 2 or 3\\
   \hline
   $[{\rm III}-{\rm III}-(d-2)/4]$ &4& 2& 1& $[{\rm III}^*-{\rm III}^*-(d-6)/4]$&4&2&3\\
   \hline
   $[{\rm I}_0-{\rm III}^*-(d-3)/4]$ &4&3& 0 or 3& $[{\rm I}^*_0-{\rm III}-(d-3)/4]$&4&3&2 or 1\\
   \hline
     $[{\rm I}_0-{\rm II}-(d-1)/6]$ &6&1&0 or 1 & $[{\rm I}^*_0-{\rm IV}^*-(d-7)/6]$&6 & 1& 3 or 4\\
     \hline
      $[{\rm II}^*-{\rm IV}-(d-7)/6]$ &6&1& 2 or 5 &$[{\rm II}^*-{\rm IV}-(d-7)/6]$& 6 & 1& 5 or 2\\
      \hline
          $[{\rm II}-{\rm IV}-(d-3)/6]$ &6& 3& 1 or 2&$[{\rm II}^*-{\rm IV}^*-(d-9)/6]$& 6& 3& 4 or 5\\
          \hline
              $[{\rm I}_0-{\rm II}^*-(d-5)/6]$ &6&5&0 or 5&$[{\rm I}^*_0-{\rm IV}-(d-5)/6]$& 6& 5& 3 or 2\\
              \hline
               $[{\rm II}-{\rm IV}^*-(d-5)/6]$ &6& 5& 1 or 4& $[{\rm II}-{\rm IV}^*-(d-5)/6]$& 6& 5& 4 or 1\\
               \hline
                 $[{\rm II}^*-{\rm III}-(d-13)/12]$ &12&1& 3 or 10&$[{\rm IV}-{\rm III}^*-(d-13)/12]$&12&1& 9 or 4\\
                 \hline
                   $[{\rm II}-{\rm III}-(d-5)/12]$ & 12 & 5& 2 or 3 & $[{\rm IV}^*-{\rm III}^*-(d-17)/12]$& 12& 5& 8 or 9\\
                   \hline
                     $[{\rm IV}-{\rm III}-(d-7)/12]$ &12 & 7& 3 or 4 & $[{\rm II}^*-{\rm III}^*-(d-19)/12]$&12&7&9 or 10\\
                     \hline
                       $[{\rm IV}^*-{\rm III}-(d-11)/12]$ & 12 & 11& 3 or 8&$[{\rm II}-{\rm III}^*-(d-11)/12]$&12&11& 9 or 2\\
                       \hline
   \end{tabular}}
   \begin{itemize}
   \item[ii)] If $2\nmid\nu(J_2)$, then the reduction type of $X^{\chi}$ is given in the following table.
   \begin{center}
\begin{tabular}{|c|c|c||c|}
  \hline
  type($X_k$) & $n$&$r$ mod $n/2$ &type($X^{\chi}_k$)\\
  \hline
   $[2{\rm I}_0-r]$ & 2&&$[2{\rm I}_0-r]$ \\
   \hline
    $[2{\rm I}^*_0-(r-1)/2]$ &4&&$[2{\rm I}^*_0-(r-1)/2]$\\
   \hline
    $[2{\rm IV}-(r-1)/3]$ &6&1&$[2{\rm IV}-(r-1)/3]$\\
   \hline
    $[2{\rm IV}^*-(r-2)/3]$ &6&2&$[2{\rm IV}^*-(r-2)/3]$\\
   \hline
    $[2{\rm III}-(r-1)/4]$ &8&1&$[2{\rm III}-(r-1)/4]$\\
   \hline
    $[2{\rm III}^*-(r-3)/4]$&8&3 &$[2{\rm III}^*-(r-3)/4]$\\
   \hline
    $[2{\rm II}-(r-1)/6]$ &12&1&$[2{\rm II}-(r-1)/6]$\\
   \hline
    $[2{\rm II}^*-(r-5)/6]$ &12&5&$[2{\rm II}^*-(r-5)/6]$\\
   \hline
\end{tabular}
\end{center}
\end{itemize}
\begin{itemize}
\vskip5pt
\item[b)] If $L/K$ is tamely ramified and $\CC_k$ is the union of two rational curves intersecting in a unique point, then the reduction type of $X^{\chi}$ is given in the following table.
\begin{center}
\begin{tabular}{|c|c|c|c||c|}
  \hline
  type($X_k$) &$n$&$\nu(J_2)$&$d$ mod $n$&type($X^{\chi}_k$)\\
  \hline
   $[{\rm I}_{d_1}-{\rm I}_{d_2}-d]$ & 1&&&$[{\rm I}^*_{d_1}-{\rm I}^*_{d_2}-(d-1)]$ \\
   \hline
$[{\rm I}^*_{d_1/2}-{\rm I}^*_{d_2/2}-(d-2)/2]$ & 2&even&0&$[{\rm I}_{d_1/2}-{\rm I}_{d_2/2}-d/2]$\\
\hline
$[{\rm I}_{e_1/2}-{\rm I}^*_{e_2/2}-(d-1)/2]$ &2&even&1&$[{\rm I}_{e_1/2}-{\rm I}^*_{e_2/2}-(d-1)/2]$ \\
\hline
$[2{\rm I}_{d_1}-d/2]$ &2&odd && $[2{\rm I}_{d_1}-d/2]$\\
\hline
$[2{\rm I}^*_{d_1/2}-(d-2)/4]$ &4& odd&& $[2{\rm I}^*_{d_1/2}-(d-2)/4]$\\
\hline
\end{tabular}
\end{center}
\vskip5pt
\item[c)] If $L/K$ is tamely ramified and $\CC_k$ is the union of an elliptic curve and a rational curve intersecting in a unique point, then the reduction type of $X^{\chi}$ is given in the following table.
    \end{itemize}
\begin{center}
{\footnotesize\begin{tabular}{|c|c|c|c||c|c|c|c|}
  \hline
  type($X_k$) &$n$& $d$ mod $n$& $r$ mod $n$ &type($X^{\chi}_k$)&$n'$&$d'$ mod $n'$& $r'$ mod $n'$\\
  \hline
   $[{\rm I}_{d_1}-{\rm I}_{0}-d]$ & 1&&&$[{\rm I}_{0}^*-{\rm I}_{d_1}^*-(d-1)]$ & 2&0&\\
   \hline
   $[{\rm I}_{0}-{\rm I}_{q}^*-(d-1)/2]$ &2&1& Remark \ref{rem:reductiontype}&$[{\rm I}_{q}-{\rm I}_{0}^*-(d-1)/2]$&2&1& \\
   \hline
   $[{\rm IV}-{\rm I}_{d_1/3}-(d-1)/3]$ &3&1&&$[{\rm II}^*-{\rm I}^*_{d_1/3}-(d-4)/3]$&6&2&  \\
   \hline
   $[{\rm IV}^*-{\rm I}_{d_1/3}-(d-2)/3]$ &3&2&&$[{\rm II}-{\rm I}^*_{d_1/3}-(d-2)/3]$&6&4& \\
   \hline
   $[{\rm III}-{\rm I}_{d_1/4}-(d-1)/4]$ &4&1&0 or 1& $[{\rm III}^*-{\rm I}^*_{d_1/4}-(d-5)/4]$ &4&1& 2 or 3 \\
   \hline
    $[{\rm III}^*-{\rm I}_{d_1/4}-(d-3)/4]$ &4&3&0 or 3&$[{\rm III}-{\rm I}^*_{d_1/4}-(d-3)/4]$&4&3&2 or 1 \\
   \hline
    $[{\rm II}-{\rm I}_{d_1/6}-(d-1)/6]$ & 6& 1& 0 or 1&$[{\rm IV}^*-{\rm I}^*_{d_1/6}-(d-7)/6]$ & 6&1& 3 or 4\\
   \hline
    $[{\rm II}^*-{\rm I}_{d_1/6}-(d-5)/6]$ &6&5&0 or 5& $[{\rm IV}-{\rm I}^*_{d_1/6}-(d-5)/6]$& 6& 5& 3 or 2 \\
   \hline
   \end{tabular}}
   \end{center}
   \end{Theorem}
   \begin{Proof}
  One knows that if $\CC_k$ consists of two elliptic curves, then so does $\CC^{\chi}_k$, see Proposition \ref{Lem:StableModelofTwist}. We assume for now that $2\mid\nu(J_2)$. Unless otherwise stated, we are assuming that $\overline{\omega}$ is non-ramified since the proof for the other subcases (b), (c), (d) of Lemma \ref{lem:n-r-qnotirreducible} is similar. Recall that $d$ is the degree of singularity of the intersection point of the irreducible components of $\CC_k$, see Proposition \ref{prop:degreeofsingularity}. Since $\CC_k$ consists of two elliptic curves, one has $d_K=\nu(J_{10}J_2^{-5})/12$. We set $u=\nu(a_0J_2)$.

   When the type of $X_k$ is $[{\rm I}_0-{\rm I}_0-d]$, one has $n=1$. In view of Lemma \ref{lem:n-r-qnotirreducible}, $n'=2$, $r'=2r+1$, and $d'=2d=0$ mod $2$. Therefore, the reduction type of $X^{\chi}$ is $[{\rm I}^*_0-{\rm I}^*_0-(d'-2)/2]$.

   When $n=2$ and $d=1$ mod $2$, i.e., the reduction type of $X$ is $[{\rm I}_0-{\rm I}^*_0-(d-1)/2]$, one has that $\displaystyle d_K\in\frac{1}{2}\Z$. Thus, $n'=2$ and $d'=d$. Therefore, the reduction type of $X^{\chi}$ is $[{\rm I}_0-{\rm I}^*_0-(d'-1)/2]$.

   When $n=3$ and $d=3d_K=0$ mod $3$ which implies that the least denominator of $d_K$ is not divisible by $3$, one has $\displaystyle u\in \frac{1}{3}\Z$ and so $n'=6$, $d'=2d=0$ mod $6$, see Lemma \ref{lem:n-r-qnotirreducible}. It follows that the reduction type of $X^{\chi}$ is $[{\rm II}-{\rm II}^*-(d'-6)/6]$. When $n=3$ and $d=1$ mod $3$, one has that $\displaystyle d_K\in\frac{1}{3}\Z$. Thus, $n'=6$, $r'=2r+3$, and $d'=2d$ is $2$ mod $6$. If $r=0$ or $1$ mod $3$, then $r'=3$ or $5$ mod $6$. In other words, the reduction type of $X^{\chi}$ is $[{\rm I}^*_0-{\rm II}^*-(d'-8)/6]$. If $r=2$ mod $3$, then $r'=1$ mod $6$ and the reduction type of $X^{\chi}$ is $[{\rm II}-{\rm II}-(d'-2)/6]$. When $n=3$ and $d=2$ mod $3$, one has $\displaystyle d_K\in\frac{1}{3}\Z$. Thus, $n'=6$ and $d'=2d=4$ mod $6$. If, moreover, $r=0$ or $2$ mod $3$, then $r'=3$ or $1$ mod $6$, and the reduction type of $X^{\chi}$ is then $[{\rm I}_0^*-{\rm II}-(d'-4)/6]$, whereas if $r=1$ mod $3$, then $r'=5$ mod $6$ and the reduction type of $X^{\chi}$ is then $[{\rm II}^*-{\rm II}^*-(d'-10)/6]$.

   When $n=4$ and $d=0$ mod $4$, one has that the denominator of $d_K$ is $1$ which implies that when $\overline{\omega}$ is non-ramified, the reduction type of $X$ cannot be $[{\rm III}-{\rm III}^*-(d-4)/4]$ since $4$ does not divide the denominator of $u/6$. Therefore, we assume that $\overline{\omega}$ is regular and ramified and $f^{-1}(\overline{\omega})$ is regular. Now, according to Proposition \ref{prop:notirreducible} the least denominator of $\nu(A_5^2J_2)/8$ is $4$, which yields that $n'=4$, $d'=d$ and $r'=r+6$, see Lemma \ref{lem:n-r-qnotirreducible}, hence the reduction type of $X^{\chi}$ is the same as this of $X$, namely $[{\rm III}-{\rm III}^*-(d-4)/4]$. The same argument holds if $\overline{\omega}$ satisfies the hypotheses of (c) or (d) in Lemma \ref{lem:n-r-qnotirreducible}. Now, we assume again that $\overline{\omega}$ is non-ramified. If $d\ne 0\mod 4$, then the denominator of $d_K$ is either $2$ or $4$. But since $n=4$ is the least common denominator of $d_K$ and $u/6$, it follows that $\displaystyle d_K\in\frac{1}{4}\Z$. This yields that $n'=4$, $d'=d$ and $r'=r+6$. Hence, the reduction type of $X^{\chi}$ follows.

    When $n=6$, we are left with dealing with the cases $d=1,3$, or $5$ mod $6$. If $d=1$ or $5$ mod $6$, then $\displaystyle d_K\in\frac{1}{6}\Z$. This yields that $n'=6$, $d'=d$ and $r'=r+3$. When $d=3$ mod $6$, one sees that $2$ is the least denominator of $d_K$, i.e., $\displaystyle d_K\in\frac{1}{2}\Z$. Therefore, $3$ is a divisor of the least denominator of $u/6$. It follows that $n'=6$, $d'=d$ and $r'=r+3$.

   When $n=12$ and $d=1,5,7,11$ mod $12$, one recalls $\displaystyle d_K=\nu(J_{10}J_2^{-5})/12$ and hence $12$ is the least denominator of $d_K$. This implies that $n'=12$, $d'=d$ and $r'=r+6$.

   Now, we assume that $2\nmid \nu(J_2)$. According to Lemma \ref{lem:n-r-qnotirreducible}, one has $[L':K]=[L:K]$, $r'=r$, and $d'=2r'=d$. It follows that the reduction type of $X^{\chi}$ is the same as the reduction type of $X$.

   The proofs of b) and c) are similar to the proof of a). The interested reader may work out the details.
   \end{Proof}

For the reduction types $[2{\rm I}_{d_1}-d/2]$ and $[2{\rm I}^*_{d_1/2}-(d-2)/4]$ in Theorem \ref{thm:twoellipticcurves} b), the degrees of singularity of the two points that are not intersection points are both equal to $d_1$. In Theorem \ref{thm:twoellipticcurves} c), since one irreducible component of $\CC_k$ is singular whereas the other component is smooth, it follows that each of the irreducible components is globally fixed under the action of $\Gal(L/K)$. This implies that $2\mid \nu(J_2)$, see \cite[Proposition 4.3.2]{Liumodelesminimaux}. Moreover, $d$ is the degree of singularity of the intersection point while $d_1$ is the degree of singularity on the rational curve.

\begin{Remark}\label{rem:reductiontype}
 Let the stable model of $C$ be the union of two irreducible components $E_1$ and $E_2$ where $\overline{\omega}\in f(E_1)$. Assume that one component is a rational curve and the other component is an elliptic curve such that $[L:K]=2$ and $d=1$ mod $2$. The reduction type of the minimal regular model $X$ of $C$ is either $[{\rm I}_{0}-{\rm I}_{q}^*-(d-1)/2]$, where $q=\nu(J_2J_{10}I_{12}^{-1})$, if $E_1$ is smooth and $r$ is even or $E_1$ is singular and $r$ is odd. Otherwise, the reduction type of $X$ is $[{\rm I}_{q}-{\rm I}_{0}^*-(d-1)/2]$, see \cite[Remarque 4.4]{Liumodelesminimaux}.
\end{Remark}

\subsection{$L/K$ is wildly ramified}
In this section we assume that $C$ admits its stable model after a wildly ramified base extension. We determine the minimal regular model of a quadratic twist of $C$ by a character $\chi$ associated with the quadratic field $K(\sqrt{D})$, $\nu(D)=1$.

The following propositions can be found in \cite[\S 5]{Liumodelesminimaux}.
 \begin{Proposition}
 \label{prop:char3}
   Suppose $3\mid[L:K]$. One has the following properties:
   \begin{itemize}
   \item[a)] $C$ is described by an equation of the form $z^2=a_0Q(u)$ with $Q(u)=(u^3+c_1u^2+c_2u+c_3)^2+c_4u^2+c_5u+c_6\in R[u]$, with $\nu(c_3)=1$ or $2$;
   \item[b)] We set $N=\min\{3\nu(c_i)-i\nu(c_3)|4\le i\le 6\}$. The reduction type of $X$ is $[{\rm III}_N]$ if $2\mid \nu(a_0)$ or $[{\rm III}^*_N]$ otherwise.
   \end{itemize}
   \end{Proposition}
\begin{Proposition}
\label{prop:char5}
   Suppose $5\mid[L:K]$. One has the following properties:
   \begin{itemize}
   \item[a)] $C$ is described by an equation of the form $z^2=b_0u^6+b_1u^5+b_2u^4+\ldots+b_6\in R[u]$, with $b_0\in\mathfrak{m}, b_1\in R^{\times}, 1\le \nu(b_6)\le 9$ and $\nu(b_6)\ne 5$;
   \item[b)] If $\nu(b_6)=2m$, the reduction type of $X$ is $[{\rm IX}-m]$. If $\nu(b_6)=2m-1$, the reduction type is $[{\rm VIII}-m]$ if $m\le 2$, or $[{\rm VIII}-(m-1)]$ if $m\ge 4$.
   \end{itemize}
   \end{Proposition}

\begin{Corollary}
\label{cor:widelyramified}
Let $C$ be a hyperelliptic curve defined over $K$. Let $K(\sqrt{D})/K$ be a quadratic extension whose associated character is $\chi$, and $\nu(D)=1$. Let $X$ and $X^{\chi}$ be the minimal regular models of the curves $C$ and its quadratic twist by $\chi$, $C^{\chi}$.
\begin{itemize}
\item[(a)] If $\Char k=3$ and the reduction type of $X$ is either $[{\rm III}_N]$ or $[{\rm III}^*_N]$, then the reduction type of $X^{\chi}$ is $[{\rm III}^*_N]$ or $[{\rm III}_N]$, respectively.
    \item[(b)] If $\Char k=5$ and the reduction type of $X^{\chi}$ is given in the following table:
\end{itemize}
{\footnotesize\begin{center}
    \begin{tabular}{|c|c|c|c|c|c|c|c|c|}
      \hline
      type($X_k$) & $[{\rm IX}-1]$&$[{\rm IX}-2]$& $[{\rm IX}-3]$& $[{\rm IX}-4]$& $[{\rm VIII}-1]$&$[{\rm VIII}-2]$&$[{\rm VIII}-3]$ &$[{\rm VIII}-4]$ \\
      \hline
      type($X^{\chi}_k$)& $[{\rm VIII}-3]$&$[{\rm VIII}-4]$&$[{\rm VIII}-1]$ &$[{\rm VIII}-2]$&$[{\rm IX}-3]$&$[{\rm IX}-4]$& $[{\rm IX}-1]$& $[{\rm IX}-2]$ \\
      \hline
    \end{tabular}
    \end{center}}
\end{Corollary}
\begin{Proof}
(a) follows from the fact that if is $C$ is defined by $z^2=a_0Q(u)$ with $Q(u)=(u^3+c_1u^2+c_2u+c_3)^2+c_4u^2+c_5u+c_6\in R[u]$, with $\nu(c_3)=1$ or $2$, then $C^{\chi}$ is defined by $z^2=Da_0Q(u)$. Now, since $\nu(a_0D)=1+\nu(a_0)$, the result follows using Proposition \ref{prop:char3} (b).

(b) Proposition \ref{prop:char5} indicated that $C$ is defined by an equation of the form $z^2=P(x)=b_0u^6+b_1u^5+b_2u^4+\ldots+b_6\in R[u]$, with $b_0\in\mathfrak{m}, b_1\in R^{\times}, 1\le \nu(b_6)\le 9$ and $\nu(b_6)\ne 5$. Since $C^{\chi}$ is defined by $z^2=DP(u)$, one may replace $u$ with $t^{-1}u$ and $z$ with $t^{-2} z$ to obtain the equation $t^{-4}z^2=D(b_0t^{-6}u^6+b_1t^{-5}u^5+b_2t^{-4}u^4+\ldots+b_6)$; or equivalently $z^2=P'(u)=b_0'u^6+b_1'u^5+b_2'u^4+\ldots+b'_6t^4$ where
\[b_0'=Db_0/t^2,b_1'=Db_1/t,b_2'=Db_2,b_3'=Dtb_3,b_4'=Dt^2b_4,b_5'=Dt^3b_5,b_6'=Dt^4b_6.\]
Therefore, $P'(u)\in R[u]$ since $b_0\in\mathfrak{m}$ and $\nu(D)=1$. Furthermore, $b_1'\in R^{\times}$. In fact, one may, and will, assume that $\nu(b_0')>0$, since otherwise one replaces $u$ with $\displaystyle\frac{u}{1-(b_0/b_1)u}$ and $z$ with $\displaystyle \frac{z}{\left(1-(b_0/b_1)u\right)^3}$ in order to obtain the equation $z^2=b_0''u^6+b_1''u^5+b_2''u^4+\ldots+b''_6t^4$ where $\nu(b_1'')=\nu(b_1')=0$, $\nu(b_i)>0$ for $i\ne 1$, and $b_6''=b_6'=Dt^4b_6$.

If $\nu(b_6)=1,2,3,4$, then $\nu(b_6')=6,7,8,9$, respectively. In other words, if the reduction type of $X$ is either $[{\rm VIII}-1]$, $[{\rm IX}-1]$, $[{\rm VIII}-2]$, or $[{\rm IX}-2]$, then the reduction type of $X^{\chi}$ is $[{\rm IX}-3]$, $[{\rm VIII}-3]$, $[{\rm IX}-4]$, or $[{\rm VIII}-4]$, respectively.

If $\nu(b_6)=6,7,8,9$, then $\nu(b_6')=11,12,13,14$, respectively. Since $\nu(b_6')>9$, one then may replace $u$ with $t^2u$ and $z$ with $t^5z$. Therefore, $\nu(b_6')=1,2,3,4$, respectively. In particular, if the reduction type of $X$ is either $[{\rm IX}-3]$, $[{\rm VIII}-3]$, $[{\rm IX}-4]$, or $[{\rm VIII}-4]$, then the reduction type of $X^{\chi}$ is
$[{\rm VIII}-1]$, $[{\rm IX}-1]$, $[{\rm VIII}-2]$, or $[{\rm IX}-2]$, respectively.
\end{Proof}
\subsection*{Acknowledgements}
 The author would like to thank Dino Lorenzini for reading an earlier draft of the paper, and several comments that improved the manuscript. We are also grateful to the referee for their thorough reading of the manuscript and for several comments and suggestions.

\end{document}